\documentclass[12pt]{article}

\usepackage{amssymb}
\usepackage{amsthm}

\usepackage{amsmath}

\begin{document}

\newtheorem{theorem}{Theorem}
\newtheorem{lemma}{Lemma}
\newtheorem{proposition}{Proposition}
\newtheorem{corollary}{Corollary}
\theoremstyle{remark}
\newtheorem{remark}{Remark}

\title{On the Convergence of Bayesian Regression Models}


\author{Yuao Hu \\Division of Mathematical Sciences, SPMS\\
        Nanyang Technological University\\
        Singapore 637371}
\maketitle
\begin{abstract}
We consider heteroscedastic nonparametric regression models, when both the mean function and variance function are unknown and to be estimated with nonparametric approaches. We derive convergence rates of posterior distributions for this model with different priors, including splines and Gaussian process priors. The results are based on the general ones  on the rates of convergence of posterior distributions for independent, non-identically distributed observations, and are established for both of the cases with random covariates, and deterministic covariates. We also illustrate that the results can be achieved for all levels of regularity, which means they are adaptive.
\end{abstract}


\section{Introduction}
\label{Introduction}
The posterior distribution is said to be consistent if the posterior probability of any small neighborhood of the true parameter value converges to one.  In recent years, many results, giving condition, under which the posterior distribution is consistent have appeared, especially under the situation that the parameter spaces are in finite-dimensional. For example, Barron et al \cite{Barron99} gave necessary and sufficient conditions for the posterior consistency , and results were then specialized to weak and $L_{1}$ neighborhoods from Kullback-Leibler neighborhoods. For details, we refer the reader to \cite{Barron99}; \cite{Ghosal99};  The consistency of posterior distributions in nonparametric Bayesian inference has received quite a lot of attention ever since 1986, when Diaconis and Freedman gave counterexample to argue that Bayesian methods sometimes can not work. On the positive side, consistency has been demonstrated on many models \cite{walker03,walker07, ga06,lian07a,lian07b,lian09rates,lian10,wu08kullback,ghosal08ejs}

In nonparametric Bayesian analysis, we have an independent sample $Y_{1},\cdots, Y_{n}$ from 
a distribution $P_{0}$ with density $p_{0}$ with respect to some measure on the sample space $(\mathcal{Y, B})$. The model space is denoted by $\mathcal{P}$ which is known to contain the true distribution $P_{0}$. Given some prior distribution $\Pi$ on $\mathcal{P}$, the posterior is a random measure given by

$$\Pi_{n}(A|Y_{1},\cdots, Y_{n})=\frac{\int_{A}\Pi_{i=1}^{n}p(Y_{i})d\Pi(P)}{\int\Pi_{i=1}^{n}p(Y_{i})d\Pi(P)}.$$
For ease of notation, we will omit the explicit conditioning and write $\Pi(A)$ for 
the posterior distribution. We say that the posterior is consistent if
$$\Pi_{n}(P: d(P,P_{0})>\epsilon)\to 0\  in\ P_{0}^{n}\ probability,$$
for any $\epsilon>0$, where $d$ is some suitable distance function between probability measures.

Furthermore, issues of rates of convergence are of interests on.  We say 
the rate is at least $\epsilon_{n}$ if for a sufficiently large constant M
$$\Pi_{n}(P: d(P,P_{0})>M\epsilon_{n})\to 0\  in\ P_{0}^{n}\ probability,$$
where $\epsilon_{n}$ is a positive sequence decreasing to zero. Ghosal and van der Vaart \cite{Ghosal00}; presented general results on the rates of convergence of the posterior measure , and \cite{Ghosal07}  then generalized the results to case even the observations are not i.i.d, which is useful for the model considered in this article. 

For Bayesian nonparametric regression models, one of the common approaches is through the splines basis expansion for regression functions, Ghosal and van der Vaart \cite{Ghosal00} gave the posterior consistency rate for regression model with unknown mean function and normal distributed error variable with zero means and known variances $\sigma^{2}$, using this approach. T.Choi and M. Schervish\cite{ConReg07} provided sufficient conditions for posterior consistency in nonparametric regression problems with homogenous Gaussian errors with unknown level by constructing tests that separate from the outside of the suitable neighborhoods of the parameter. Amewou-Atisso, Ghosal, Ghosh and Ramamoorthi \cite{Consemi03} presented a posterior consistency analysis for linear regression problems with an unknown error distribution which is symmetric about zero. Besides, both papers did not consider the rates of convergence. 

In this paper, we give the convergence rates for heteroscedastic nonparametric regression models, when we use nonparametric methods to estimate the unknown variance function and the unknown mean function simultaneously. Besides, as in \cite{Ghosal99}, we also deal with two types of covariates either randomly sampled from a probability distribution or fixed in advance. When the covariate values in one-dimensional, we use the approach of splines basis expansion for regression functions and give the convergence rate. For high dimensional cases, we use rescaled smooth Gaussian fields as priors for multidimensional functions to get the result.

Using Gaussian process in the context of density estimation is another common approach in Bayesian nonparametric analysis. It is first used by Leonard\cite{Leonard78} and Lenk \cite{lenk88}. Recently, many results on posterior consistency are induced by the Gaussian process prior, such as in  \cite{Leonard78}, and  \cite{ConGau08}. Van der Vaart and van Zanten \cite{RepGau08} derived the rates of contraction of posterior distributions on nonparametric or semiparametric models based on Gaussian processes and showed that the rates depend on the position of the true parameter associated with the reproducing kernel Hilbert space of the Gaussian process and the small ball probabilities of the Gaussian process. With rescaled smooth Gaussian fields as priors, they\cite{AdaGau09} extended the results to be fully adaptive to the smoothness.

The rest of the paper is organized as follows. In Section 2, we describe the regression model. In Section 3, we give the main results of the posterior convergence rates with splines basis expansion approach and Gaussian process approach. Section 4 contains proofs with some lemma left to the Appendix. We discuss about the results and some directions on future work in section 5.


\section{The model}
\label{model}
We consider the heteoscedastic nonparametric regression model, where a random response $y$ corresponding to a covariate vector $\mathbf{x}$ taking values in a compact set $T\subset R^{d}$, without loss of generality, we assume that $T=[0,1]^{d}.$

To be specific, the regression model we consider here, is the following:
\begin{equation}
\label{model}
\begin{split}
&y_{i}=\eta(\mathbf{x_{i}})+V^{1/2}(\mathbf{x_{i}})\epsilon_{i},\\
&\epsilon_{i}\sim N(0,1),\\
&\eta(.)\sim\Pi_{1},\\
&f(.)=\log V(.)\sim\Pi_{2},\\
\end{split}
\end{equation}

when $\eta(.)$ is the mean function, and $V(.)$ is the variance function.  Let $\Theta^{(g)}$ be the abstract measure space where the function $g(x)$ ($g(x)$ indicates $\eta(x)$, or $ f(x)$) belongs to, with respect to a common $\sigma$-finite measure. With the assumption $\Theta^{(\eta)}$ and $\Theta^{(f)}$ are independent, we define the jointly parameter space $\Theta$ the product space of  $\Theta^{(\eta)}$ and  $\Theta^{(f)}$. Since the parameter space is infinite-dimensional, we consider a sieve $\Theta_{n}$ growing eventually to the space of $\Theta$,with $\Theta_{1}\subseteq\cdots, \subseteq \Theta_{n}\subseteq \Theta$ and $\cup\Theta_{n}=\Theta$. We model the unknown function $\eta(.)$ and $f(.)$ with suitable prior distributions $\Pi_{n}^{(\eta)}$ and $\Pi_{n}^{(f)}$ on their parameter sieve spaces, respectively.

\section{The main results }
\label{main result}
In this section, we give the rates of convergence of the nonparametric regression model described in Section~\ref{model}, The parameter is $\theta=(\eta, V$) with $\theta_{0}=(\eta_{0}, V_{0})$ being the true functions.

Let  $\mathcal{P}_{\eta,V}$ be the distribution of $y$. To be specific, for our model 
\begin{equation}
\label{distribution}
\begin{split}
\mathcal{P}_{\eta,V}(y|x)=\frac{1}{\sqrt{2\pi V(\mathbf{x})}}exp(-\frac{(y-\eta(\mathbf{x}))^{2}}{2V(\mathbf{x})}).
\end{split}
\end{equation}

We use $d_{n}^{2}$ to denote the squares of the Hellinger distances. It means, for  random covariates and fixed covariates 
\begin{equation}
\label{distribution}
\begin{split}
d^{2}_{n}(\mathcal{P}_{\eta_{1},V_{1}},\mathcal{P}_{\eta_{2},V_{2}})=\int\int(\mathcal{P}_{\eta_{1},V_{1}}^{\frac{1}{2}}-\mathcal{P}_{\eta_{2},V_{2}}^{\frac{1}{2}})^{2}\ dy\ dQ(\mathbf{x}).
\end{split}
\end{equation}
For random covariates, $Q(\mathbf{x})$ denotes the distribution function of $\mathbf{x}$, and for fixed covariates, it is the empirical probability measure of the design points, which is defined by$P_{n}^{\mathbf{x}}=n^{-1}\sum^{n}_{i=1}\delta_{\mathbf{x_{i}}}$.

The Kullback-Leibler divergence and variance divergence of $ P_{\eta_{1},V_{1}}$ and $P_{\eta_{2},V_{2}}$ for fixed $\mathbf{x}$ are defined in the following way:
\begin{equation}
\label{KL}
\begin{split}
&K_{\mathbf{x}}(\mathcal{P}_{\eta_{1},V_{1}},\mathcal{P}_{\eta_{2},V_{2}})=\int \mathcal{P}_{\eta_{1},V_{1}}\ \log (\frac{\mathcal{P}_{\eta_{1},V_{1}}}{\mathcal{P}_{\eta_{2},V_{2}}})\ dy;\\
&Var_{\mathbf{x}}(\mathcal{P}_{\eta_{1},V_{1}},\mathcal{P}_{\eta_{2},V_{2}})=\int\mathcal{P}_{\eta_{1},V_{1}}\ (\log (\frac{\mathcal{P}_{\eta_{1},V_{1}}}{\mathcal{P}_{\eta_{2},V_{2}}})-K_{\mathbf{x}}(\mathcal{P}_{\eta_{1},V_{1}},\mathcal{P}_{\eta_{2},V_{2}}))^{2}\ dy.
\end{split}
\end{equation}\\
For the specific model in section 2:
\begin{equation}
\label{KLspecial}
\begin{split}
&K_{\mathbf{x}}(\mathcal{P}_{\eta_{1},V_{1}},\mathcal{P}_{\eta_{2},V_{2}})=\frac{1}{2}\log\frac{V_{2}(\mathbf{x})}{V_{1}(\mathbf{x})}-\frac{1}{2}(1-\frac{V_{1}(\mathbf{x})}{V_{2}(\mathbf{x})})+\frac{1}{2}\frac{[\eta_{1}(\mathbf{x})-\eta_{2}(\mathbf{x})]^{2}}{V_{2}(\mathbf{x})};\\
&Var_{\mathbf{x}}(\mathcal{P}_{\eta_{1},V_{1}},\mathcal{P}_{\eta_{2},V_{2}})=2[-\frac{1}{2}+\frac{1}{2}\frac{V_{1}(\mathbf{x})}{V_{2}(\mathbf{x})}]^{2}+[\frac{V_{1}(\mathbf{x})}{V_{2}(\mathbf{x}))}[\eta_{1}(\mathbf{x})-\eta_{2}(\mathbf{x})]]^{2}.
\end{split}
\end{equation}\\
Correspondently, the average Kullback-Leibler divergence and variance divergence are in the forms of 
\begin{equation}
\label{KLaverage}
\begin{split}
&K(\mathcal{P}_{\eta_{1},V_{1}},\mathcal{P}_{\eta_{2},V_{2}})=\int K_{\mathbf{x}}(\mathcal{P}_{\eta_{1},V_{1}},\mathcal{P}_{\eta_{2},V_{2}})\ dQ(\mathbf{x});\\
&Var(\mathcal{P}_{\eta_{1},V_{1}},\mathcal{P}_{\eta_{2},V_{2}})=\int Var_{\mathbf{x}}(\mathcal{P}_{\eta_{1},V_{1}},\mathcal{P}_{\eta_{2},V_{2}})\ dQ(\mathbf{x}).
\end{split}
\end{equation}\\
In the remainder of the article, let $||.||_{n}$ stand for the norm on $L_{2}(Q)$, $||.||_{\infty}$ denotes the supreme norm .
\subsection{ Splines}\label{Splines}
In this section we give the convergence rates to prior distributions on spline models for regression  functions. We restrict ourselves to the one-dimensional case here, though for higher dimensions case, tensor splines can be used.

The basic assumption for the true densities of the mean function and variance function is that they belong to the  H$\ddot{o}$lder spaces $C^{\alpha}[0,1]$ and $C^{\gamma}[0,1]$, respectively, where $\alpha, \gamma >0$ could be fractional. The H$\ddot{o}$lder space $C^{\alpha}[0,1]$ is constructed by all functions that have $\alpha_{0}$ derivatives, with $\alpha_{0}$ being the greatest integer less than $\alpha$ and $\alpha_{0}$th derivative being Lipschitz of order $\alpha-\alpha_{0}$.

Throughout this article, we fix an order q, which is a natural number satisfied $q\geq max\{\alpha, \gamma\}$. A B-spline basis function of order $q$ consists of $q$ polynomial pieces of degree $q-1$, which are $q-2$ times continuously differentiable throughout [0,1]. To approximate a function on [0,1], we partition the interval [0,1] into $K_{n}$ subintervals $((k-1)/K_{n}, k/K_{n}]$ for $k=1,2,\cdots, K_{n}$, with $\{K_{n}\}$ being a sequence of natural numbers increasing to infinity as n goes to infinity. Each subinterval $((k-1)/K_{n}, k/K_{n}]$ is approximated by a polynomials of degree strictly less than $q$. The number of basis functions needed is $J_{n}=(q+K_{n}-1)$. The basis functions can be denoted as $B_{j}$, with $j=1,2,\cdots J_{n}$. Thus, the space of splines of order $q$ is a  $J_{n}$-dimensional linear space, consisted by all functions from $[0,1]$ to $\mathbb{R}$ in form of $g=\sum_{j=1}^{J^{n}}\beta_{j}B_{j}$. As in \cite{Ghosal07}, the B-splines satisfy (i) $B_{j} \geq 0$, $j=1,2,\cdots J_{n}$, (ii) $\sum_{j=1}^{J_{n}}B_{j}=1$, (iii) $B_{j}$ is supported inside an interval of length $q/K_{n}$ and (iv) at most q of $B_{1}, B_{2}, \cdots, B_{J_{n}}$ are nonzero at any given x.

We denote g= $f $ or $\eta$, and put prior on g by a prior on $\mathbb{\beta}=(\beta_{1}, \cdots, \beta_{J_{n}})^{T}$, the spline coefficients, where $g$ is represented as $g_{\beta}(x)=\beta^{T}B(x)$. Let $\Pi_{n}^{(g)}$ be priors induced by a multivariate normal distribution $N_{J_{n}}(0, I)$ on the spline coefficients.

We also assume the regressors are sufficiently regularly distributed, by satisfying the condition expressed in the following term 
\begin{equation}
\label{regu}
J_{n}^{-1}||\beta||^{2} \lesssim \beta^{T} \Sigma \beta \lesssim J_{n}^{-1}||\beta||^{2},
\end{equation}
 where $\Sigma=(\int B_{i}B_{j}\ d\ Q)$, $||.||$ is the Euclidean norm on $\mathbb{R}^{J_{n}}$.\\
\begin{theorem}
\label{theorem1}
 Assume that $\eta_{0}\in C^{\alpha}[0,1], V_{0}\in C^{\gamma}[0,1]$ for some $\alpha, \gamma \geq \frac{1}{2}$, $V_{0}$ is away from 0, and (\ref{regu}) holds. Let $\Pi_{n}^{(\eta)}$ and $\Pi_{n}^{(f)}$ be  priors of $\eta$ and $f$ both induced by $N_{J_{n}}(0, I)$ on the spline coefficients. If $$J_{n} \sim \min\{(n/\log n)^{1/(1+2\alpha)}, n^{1/(2+2\gamma)}\},$$then the posterior converges at the rate $$ \epsilon_{n}\sim \max\{ (n/\log n)^{-\alpha/(1+2\alpha)}, n^{-\gamma/(2+2\gamma)}\},$$
 relative to $d_{n}$.
\end{theorem}

Usually, we can view $J_{n}$ to be a sequence of random variables with a prior distributions. It can be prove that the posterior can convergence at the same rate.\\
\begin{corollary}
\label{corollary1}
Assume that $\eta_{0}\in C^{\alpha}[0,1], V_{0}\in C^{\gamma}[0,1]$ for some $\alpha, \gamma \geq \frac{1}{2}$, $V_{0}$ is away from 0, and (\ref{regu}) holds. Let $\Pi_{n}^{(\eta)}$ and $\Pi_{n}^{(f)}$ be  priors of $\eta$ and $f$ both induced by $N_{J_{n}}(0, I)$ on the spline coefficients. $J_{n}$ is a sequence of geometric distributed random variables with successful probability $p_{n}$ satisfying $ p_{n}^{k_{n}-1}(1-p_{n})=e^{-n\epsilon^{2}_{n}}$, with $k_{n}=\lfloor\min\{( n/\log n)^{1/(2\alpha+1)}, n^{1/(2+2\gamma)}\}\rfloor$, $ \epsilon_{n}\sim \max\{ (n/\log n)^{-\alpha/(1+2\alpha)}, n^{-\gamma/(2+2\gamma)}\}. Then, $ the posterior convergence rate is $\epsilon_{n}$, relative to $d_{n}$.
\end{corollary}

\subsection{Gaussian process prior}\label{Gaussian process prior}

For higher dimensional case, we employ prior distributions, constructed by rescaling smooth Gaussian random field. Let $\Theta$ be $C[0,1]^{d}$, the space of all continuous functions defined on $[0,1]^{d}$. As in \cite{AdaGau09}, we set $W^{(g)}=(W_{\mathbf{x}}^{(g)}: \mathbf{x}\in\mathbb{R}^{d})$ to be a centered, homogeneous Gaussian random field with covariance function of the form, for a given continuous function $\phi$:
$$EW_{s}^{(g)}W_{t}^{(g)}=\phi(s-t).$$To be specific, we choose  $W^{(g)}=(W_{x}^{(g)}: \mathbf{x}\in\mathbb{R}^{d})$ to be the squared exponential process, which is the centered Gaussian process with covariance function 
$$EW_{s}^{(g)}W_{t}^{(g)}=exp(-||s-t)||^{2}),$$ where $||.||$ is the Euclidean norm on $\mathbb{R}^{d}$.

Let $A$ be a random variable defined on the same probability space as $W^{(g)}$ and independent of $W^{(g)}$. Here we assume $A^{d}$ possesses a Gamma distribution. $W^{(g)A}$ is used to denote the rescaled process $x\to W_{Ax}$ restricted on $[0,1]^{d}$, which can be considered as a Borel measurable map in the space $C[0,1]^{d}$, with the uniform norm $||.||_{\infty}$, as showed in \cite{AdaGau09}.\\
\begin{theorem}
\label{theorem2}
Assume that $\eta_{0}\in C^{\alpha}[0,1]^{d}, V_{0}\in C^{\gamma}[0,1]^{d}$ for some $\alpha, \gamma \geq \frac{1}{2}$, $V_{0}$ is away from 0. We consider the prior on g is (g denotes f or $\eta$) $W^{(g)A}$, which is the restricted and rescaled squared exponential process with $A^{d}$ a Gamma distributed random variable. Then, the posterior converges at the rate $$\epsilon_{n}=\max\{ n^{-\alpha/(d+2\alpha)}(\log n)^{(d+1)\alpha/(2\alpha+d)}, n^{-\gamma/(d+2\gamma)}(\log n)^{(d+1)\gamma/(2\gamma+d)}\}$$ relative to $d_{n}$.
\end{theorem}

The proof can be found in Section 4. Also, this rate of contraction is not minimax. By choosing a different prior for $A$, the power $(d+1)\alpha/(2\alpha+d)$ of the logarithmic factor can be improved. Though the prior does not depend on $\alpha$ and $\gamma$, the convergence rate is true for any level of $\alpha$, and $\gamma$. In this sense, it is rate-adaptive. 

If we do not consider about the property of adaption or the regularity levels are known, we can find the minimax rate by using proper priors. \\
\begin{corollary}
\label{corollary2}
Assume that $\eta_{0}\in C^{\alpha}[0,1], V_{0}\in C^{\gamma}[0,1]$ for some $\alpha, \gamma \geq \frac{1}{2}$, $V_{0}$ is away from 0. For simplicity, we only consider the one-dimensional situation for simplicity. We denote $W^{(g)}$ to be a standard Brownian motion and $Z_{0},\cdots Z_{k_{g}}$ independent standard normal random variables. We consider the prior on g is the process $x\to I_{0+}^{k_{g}}W^{(g)}_{x}+\sum_{i=1}^{k_{g}}Z_{i}x^{i}/i!$, where $I_{0+} W$ denotes $x\to \int_{0}^{x}W(x)dx$, and $I_{0+}^{k}W$ denotes $ I_{0+}^{1}(I_{0+}^{k-1}W)$. Then, the posterior converges at the rate $$\max\{n^{-\alpha/(2k_{\eta}+2)}, n^{-\gamma/2k_{f}+2}\}$$When $\gamma=k_{f}+1/2$ and $\alpha=k_{\eta}+1/2$, $$\epsilon_{n}=\max\{ n^{-\alpha/(1+2\alpha)}, n^{-\gamma/(1+2\gamma)}\}$$
which is the minimax rate. 
\end{corollary}

This example shows, for the case $\alpha$, and $\gamma$ are known, we can use the above specific Gaussian process prior to get the minimax rate. However, this is not optimal for all level of $\alpha$ and $\gamma$, so other choice of $k_{g}$ will corresponds to under-or-over-smoothed prior. 

\section{ The proofs for the main results }\label{proofs}

In preparation for the proofs of the main results, we first collect some lemmas, which are used to bound the average hellinger distance entropy, Kullback-Leibler divergence and variance divergence with the $L_{2}$ norm of the regression functions.\\
\begin{lemma}
\label{LemmaEn}
The average hellinger distance entropy of the product space $\Theta_{n}$ can be bounded by a multiple of the summation of $||.||_{n}$-entropy of $\Theta_{n}^{(\eta)}$ and $\Theta_{n}^{(f)}$, reminding that $f=\log V$, which means 
\begin{equation}
\label{entropy}
\log N(3\epsilon, \Theta_{n}, d_{n}) \lesssim\ \log N(\epsilon/e^{N_{n}}, \Theta_{n}^{(\eta)},  ||.||_{n})+\log N(\epsilon, \Theta_{n}^{(f)},  ||.||_{n}).
\end{equation}
\end{lemma}

With this lemma, the $\epsilon$-covering number relative to $d_{n}$-metric can be estimated that with relative to $L_{2}$-metric. \\\\
\begin{lemma}
\label{LemmaKL}
Under the assumption that both $f_{1}$ and $f_{2}$ are uniformly bounded by a constant N, 
\begin{equation}
\label{KLestimation}
\begin{split}
&K(\mathcal{P}_{\eta_{1},V_{1}},\mathcal{P}_{\eta_{2},V_{2}})\leq (1+ e^{2N})(||\eta_{1}-\eta_{2}||_{n}^{2}+||(f_{1}-f_{2})||_{n}^{2});\\
&Var(\mathcal{P}_{\eta_{1},V_{1}},\mathcal{P}_{\eta_{2},V_{2}}) \leq e^{4N}(||\eta-\eta_{0}||_{n}^{2}+||f-f_{0}||_{n}^{2}).
\end{split}
\end{equation}
\end{lemma}

We use this lemma to estimate the prior concentration probability. The proofs can be found in the Appendix.\\

\subsection{Proof for theorem 1 }\label{proof No.1}
We consider sieve $\Theta_{n}=\Theta_{n}^{(f)}\times\Theta_{n}^{(\eta)}$ where  $$\Theta_{n}^{(f)}=\{f_{\beta}\in supp\{\Pi_{n}^{(f)}\}, ||f_{\beta}||\leq N_{n}\};$$ $$\Theta_{n}^{(\eta)}=\{\eta_{\beta}\in supp\{\Pi_{n}^{(\eta)}\}, ||\eta_{\beta}||\leq M_{n}\},$$ where $supp\{\Pi_{n}\}$ means the support of $\Pi_{n}$ and $M_{n}$, $N_{n}$ are sequence of real numbers goes to infinity as n goes to infinity. Since we suppose $\eta_{0}\in C^{\alpha}[0,1], V_{0}\in C^{\gamma}[0,1]$, and $V_{0}$ is away from 0, we have $\log V_{0} \in C^{\gamma}[0,1]$, too. By the Lemma 4.1 in \cite{Ghosal00}, there exists some $\beta_{\eta_{0}}, \beta_{f_{0}}\in\mathbb{R}^{J_{n}} $(dependent on n), for the true density of $f_{0}$ and $\eta_{0}$, the basic approximation property of splines are satisfied as 
\begin{equation}
\label{approx1}
\begin{split}
&||\beta_{f_{0}}^{T}B-f_{0}||_{\infty}\leq A J_{n}^{-\gamma}||f_{0}||_{\gamma};\\
&||\beta_{\eta_{0}}^{T}B-\eta_{0}||_{\infty} \leq A^{'}J_{n}^{-\alpha}||\eta_{0}||_{\alpha},
\end{split}
\end{equation}
where $A$, and $A^{'}$ are constant.

Under the assumption of (\ref{regu}) in Theorem 1, we can use Euclidean norms on the spline coefficients to control the $L_{2}$ distance of functions, since for all $\beta, \beta^{'}\in \mathbb{R}^{J_{n}} $,
\begin{equation}
\label{approx2}
C^{-1}||\beta-\beta^{'}||\leq \sqrt{J}||g_{\beta}-g_{\beta^{'}}||_{n}\leq (C^{'})^{-1}||\beta-\beta^{'}|| 
\end{equation}
are satisfied for some constants C and $C^{'}$. 

 We verify all the conditions of general results on rates of posterior contraction (e.g. Theorem 4 of \cite{Ghosal07} ), except that the local entropy in condition (3.2) is replaced by the global entropy $ \log N(\epsilon, \Theta_{n}, d_{n})$ without affection rates. The parameter $\theta$ in Theorem 4 of \cite{Ghosal07} is ($\eta, V$) with $\theta_{0}=(\eta_{0}, V_{0})$.  
 
We start from the estimation of entropy number. We project $g_{0}$ onto the $J_{n}$-dimensional space of splines and denote the projection function $g_{\beta^{(n)}_{g}}$. Using the property of projection combined with (\ref{approx2}), we have that $\{\beta:||g_{\beta}-g_{0}||_{n}\leq \epsilon\}\subset\{\beta:||\beta-\beta^{(n)}_{g}||\leq C\sqrt{J_{n}}\epsilon\}$ for every $\epsilon>0$. For details, please refer to \cite{Ghosal07}. Thus, we can use the $C\sqrt{J_{n}}\epsilon$-covering numbers relative to Euclidean norm to bound the $\epsilon$-covering number of the set $\{\beta:||g_{\beta}-g_{0}||_{n}\}$ relative to $L_{2}$ norm. Thus, we have 
\begin{equation}
\label{EnEstimation}
N(\epsilon/3, \Theta_{n}^{(\eta)},  ||.||_{n})\lesssim N(C\sqrt{J_{n}}\epsilon, \Theta_{n}^{(\eta)},  ||.||)\lesssim (\frac{KM_{n}}{\epsilon_{n}})^{J_{n}},
\end{equation}
 where $K$ is a constant, $\eta$ can be replaced by $f$ with $M_{n}$ replaced by $N_{n}$ together. So by lemma~\ref{LemmaEn}, the entropy condition $\log N(\epsilon, \Theta_{n}, d_{n}) \lesssim n \epsilon^{2}$ is satisfied, provided $J_{n}\log M_{n}\lesssim n\epsilon_{n}^{2}$, $J_{n}N_{n}\lesssim n\epsilon_{n}^{2}$ and $J_{n}\log \epsilon_{n}^{-1}\lesssim n\epsilon_{n}^{2}$.

Then, we turn to estimate the prior concentration probability for the true density, which is in form of 
\begin{equation}
\label{concentration}
\Pi_{n}(B_{n}((\eta_{0},f_{0}),\epsilon_{n};2))=\bigg\{(\eta, V): K(\mathcal{P}_{\eta,V},\mathcal{P}_{\eta_{0},V_{0}})\leq \epsilon^{2}, Var(\mathcal{P}_{\eta,V},\mathcal{P}_{\eta_{0},V_{0}})\leq \epsilon^{2}\bigg\}.
\end{equation}
We denote $N/2$ to be $\ ||f_{0}||_{\infty}$. Under the assumption that $||f||_{\infty}\leq N$ and $\ ||f_{0}||_{\infty}\leq N$, when n is sufficiently large, 
\begin{equation}
\label{ConcentrationEs}
\begin{split}
&\Pi_{n}(B_{n}((\eta_{0},f_{0}),\epsilon_{n};2))\\
&\geq \bigg\{(\eta, V): K(P_{\eta,V},P_{\eta_{0},V_{0}})\leq \epsilon^{2}, Var(P_{\eta,V},P_{\eta_{0},V_{0}})\leq \epsilon^{2}, ||f||_{\infty}<N\bigg\}\\
&\geq \Pi_{n}(||f-f_{0}||_{n}^{2}+||\eta-\eta_{0}||_{n}^{2}\leq e^{-4N}\epsilon_{n}^{2}, ||f||_{\infty}<N)\\
&\geq\Pi_{n}^{(f)}(f:||f-f_{0}||_{n}^{2}\leq \frac{e^{-4N}}{2}\epsilon_{n}^{2},||f||_{\infty}<N)\times \Pi_{n}^{(\eta)}(\eta:||\eta-\eta_{0}||_{n}^{2}\leq 
\frac{e^{-4N}}{2}\epsilon_{n}^{2})\\
&\geq \Pr_{\beta^{T}B\in \Theta_{n}^{(f)}}(\beta: ||\beta-\beta^{(n)}_{f}||\leq e^{-2N} C^{'}\sqrt{J_{n}}\epsilon_{n},|\beta^{(n)}_{j}|<N)\\
&\times  \Pr_{\beta^{T}B\in \Theta_{n}^{(\eta)}}(\beta: ||\beta-\beta^{(n)}_{\eta}||\leq e^{-2N} C^{'}\sqrt{J_{n}}\epsilon_{n},|\beta^{(n)}_{j}|<N)\\
&\geq (\mathop{\inf}\limits_{\beta_{1}\in [-2N, 2N]} \phi(\beta_{1}))^{2}Vol(\beta: ||\beta-\beta^{(n)}_{f}||\leq e^{-2N} C^{'}\epsilon_{n})Vol(\beta: ||\beta-\beta^{(n)}_{\eta}||\leq e^{-2N} C^{'}\epsilon_{n})\\
&\gtrsim \epsilon_{n}^{2J_{n}}
\end{split}
\end{equation}
where $vol$ denotes the volume in Euclidean space and $\mathop{\inf}\limits_{\beta_{1}\in [-2N, 2N]}\phi(\beta_{1})$ represents the infimum value of density function $\phi$, which is the density function of normal distribution, constrained on the open set $[-2N, 2N]$. The second inequality is derived from lemma~\ref{LemmaKL}.$\mathop{\inf}\limits_{\beta_{1}\in (-2N, 2N)}\phi(\beta_{1})$ is a real number away from zero, which can be derived from the facts that $\phi$ is nonzero at any point belongs to $\mathbb{R}$ alone with its continuity, and  [-2N, 2N] is a compact set in $\mathbb{R}$.

To satisfy the entropy and the prior concentration conditions, it is necessary  that $J_{n}N_{n}\lesssim n\epsilon^{2}_{n}$, $J_{n}\log M_{n}\lesssim n\epsilon^{2}_{n}$, and $J_{n}\log \epsilon_{n}^{-1}\lesssim n\epsilon_{n}^{2}$ together with $\epsilon_{n}\gtrsim 2J_{n}^{-\nu}$, where $\nu=\min\{\alpha, \gamma\}$. When we set $N_{n}\sim n^{1/(2\nu+2)}$, $M_{n}\sim n$, all conditions of above are satisfied, with $$J_{n} \sim \min\{(n/\log n)^{1/(1+2\alpha)}, n^{1/(2+2\gamma)}\},$$ and $$ \epsilon_{n}\sim \max\{ (n/\log n)^{-\alpha/(1+2\alpha)}, n^{-\gamma/(2+2\gamma)}\}.$$

The left is to get the condition on which the probability assigned by prior to $\Theta_{n}$ complement is exponentially small. As we mentioned, $\eta_{\beta}=\beta^{T}B(x)$ for all $x\in[0,1]^{d}$, and $|\sum_{j}^{J_{n}}\beta_{j}B_{j}|\leq \max_{j=1}^{J_{n}}|\beta_{j}|$. Then for $t_{n}>0$, by Markov's inequality and Chernoff Bounds, we have
\begin{equation} 
\label{TailEst1}
\Pr \Big\{\mathop{\sup}\limits_{x\in[0,1]}|\sum_{j}^{J_{n}}\beta_{j}B_{j}|> M_{n}\Big\}\leq J_{n}\exp\Big(-t_{n}M_{n}+\frac{1}{2}t_{n}^{2}\Big)2\Phi(t_{n}),
\end{equation}
where $\Phi$ is the standard normal distribution function. By taking $t_{n}=M_{n}$, we have 
\begin{equation}
\label{TailEst2}
\Pr \Big\{\mathop{\sup}\limits_{x\in[0,1]}|\sum_{j}^{J_{n}}\beta_{j}B_{j}|> M_{n}\Big\}\lesssim J_{n}exp\Big(-\frac{M_{n}^{2}}{2}\Big).
\end{equation}
With the $M_{n}$, $N_{n}$, $J_{n}$ and $\epsilon_{n}$ defined as above, and n sufficiently large,
\begin{equation}
\label{TailEst3}
J_{n}\exp\Big(-\frac{M_{n}^{2}}{2}\Big)\lesssim \exp\Big(-n\epsilon_{n}^{2}\Big),
\end{equation}
and the formula replacing $M_{n}$ with $N_{n}$ are also satisfied. Thus,
\begin{equation}
\label{TailEstFinal}
\begin{split}
&\Pi_{n}(\Theta\setminus\Theta_{n})\leq \Pi_{n}^{(f)}(\Theta^{(f)}\setminus\Theta_{n}^{(f)})+\Pi_{n}^{(\eta)}(\Theta^{(\eta)}\setminus\Theta_{n}^{(\eta)})\\
&=\Pr \Big\{\mathop{\sup}\limits_{x\in[0,1]}|\sum_{j}^{J_{n}}\beta_{j}B_{j}|> M_{n}\Big\} +\Pr \Big\{\mathop{\sup}\limits_{x\in[0,1]}|\sum_{j}^{J_{n}}\beta_{j}B_{j}|> N_{n}\Big\} \\
&\lesssim \exp\Big(-n\epsilon_{n}^{2}\Big).
\end{split}
\end{equation}
The whole proof is completed.\\
\begin{remark}
When we generalize the priors of $\eta$ and $f$, which are induced by the spline coefficients, with some limitation, the convergence rate will stay unchanged. We assume the same prior $\Pi$ on each $\beta_{j}\in \mathbb{R}$, $j=1,\cdots, J_{n}$, with density function $d(\beta_{j})\in C[\mathbb{R}]$ (the set of continuous functions), which satisfies
\begin{equation}
\label{DensityCon}
\begin{split} 
&\Pi (|\beta_{j}|>M)\lesssim e^{- M^{\rho}};\\ 
& d(\beta_{j}= r )\neq 0 \ for\ any \ r\in \mathbb{R},
\end{split}
\end{equation}
where $\rho$ is a real number larger than 1. The normal distribution can be viewed as a special case satisfying (\ref{DensityCon}). Then, with $ \epsilon_{n}$, $M_{n}$ $N_{n}$, and $J_{n}$ defined as above, \ref{TailEstFinal} are not affected, since 
\begin{equation}
\label{TailGener}
\Pr \Big\{\mathop{\sup}\limits_{x\in[0,1]}|\sum_{j}^{J_{n}}\beta_{j}B_{j}|> M_{n}(N_{n}, resp.))\Big\}\lesssim J_{n}\exp\Big(-\frac{M_{n}^{\rho}}{2}\Big)\lesssim \exp\Big(-n\epsilon_{n}^{2}\Big).
\end{equation}

The prior concentration probability estimation can also be bounded  below by a multiple of the volume of a Euclidean ball. Added with the fact that priors does not affect the entropy, we finish showing that the convergence rate can keep still when we generalize the priors. 
\end{remark}

\subsection{Proof for corollary 1 }\label{proof No.2}
The proof is almost the same with that for theorem 1. We consider the sieves $\Theta_{n}=\Theta_{n}^{(f)}\times\Theta_{n}^{(\eta)}$ in the form of 
$$\Theta_{n}^{(f)}=\{f_{\beta}\in supp\{\Pi_{n}^{(f)}\}, ||f_{\beta}||\leq N_{n}, J_{n}\leq k_{n}\};$$ $$\Theta_{n}^{(\eta)}=\{\eta_{\beta}\in supp\{\Pi_{n}^{(\eta)}\}, ||\eta_{\beta}||\leq M_{n},J_{n}\leq k_{n}\},$$
where $k_{n}=\lfloor\min\{( n/\log n)^{1/(2\alpha+1)}, n^{1/(2+2\gamma)}\}\rfloor$ and $\lfloor.\rfloor$ denotes the Integral part. With (\ref{EnEstimation}), the $\epsilon_{n}$-entropy of $\Theta_{n}$ is bounded by a multiple of $(\frac{M_{n}}{\epsilon_{n}e^{-N_{n}}})^{J_{n}}\times (\frac{N_{n}}{\epsilon_{n}})^{J_{n}}$, which have been proved to be always bounded by a multiple of $e^{n\epsilon^{2}}$ with $ J_{n}\leq \lfloor\min\{( n/\log n)^{1/(2\alpha+1)}, n^{1/(2+2\gamma)}\}\rfloor$, $M_{n}\sim n$, $N_{n}\sim n^{1/(2\gamma+2)}$and $ \epsilon_{n}\sim \max\{ (n/\log n)^{-\alpha/(1+2\alpha)}, n^{-\gamma/(2+2\gamma)}\}.$

The prior concentration probability (\ref{concentration}) can be estimated in the form of 
\begin{equation*}
\begin{split}
&\Pi_{n}(B_{n}((\eta_{0},f_{0}),\epsilon; 2))\\
&=\mathop{\sum}\limits_{k=1}\limits^{k_{n}}\Pr (J_{n}=k) \Pi_{n}(B_{n}((\eta_{0},f_{0}),\epsilon;2), J_{n}=k)\\
&\geq \Pr (J_{n}=k_{n})(\mathop{\inf}\limits_{\beta_{1}\in [-2N, 2N]} \phi(\beta_{1}))^{2}Vol(\beta: ||\beta-\beta^{(n)}_{f}||\leq e^{-2N} C^{'}\epsilon)Vol(\beta: ||\beta-\beta^{(n)}_{\eta}||\leq e^{-2N} C^{'}\epsilon)\\
& \gtrsim \Pr (J_{n}=k_{n})\epsilon^{2k_{n}}\\
\end{split}
\end{equation*}
With the assumption for $p_{n}$ and the fact that we have already proved $\epsilon^{2J_{n}} \gtrsim e^{-n\epsilon^{2}}$ with $J_{n} \sim \min\{(n/\log n)^{1/(1+2\alpha)}, n^{1/(2+2\gamma)}\},$ and $ \epsilon_{n}\sim \max\{ (n/\log n)^{-\alpha/(1+2\alpha)}, n^{-\gamma/(2+2\gamma)}\}.$
we can guarantee $$p_{n}^{k_{n}-1}(1-p_{n})\epsilon_{n}^{2k_{n}}\gtrsim e^{-n\epsilon_{n}^{2}}.$$
We compute the probability of $(\Theta_{n}^{(\eta)})^{c}$ as following:
\begin{equation*}
\begin{split}
&\Pi_{n}((\Theta_{n}^{(\eta)})^{c})\\
&=\mathop{\sum}\limits_{k=1}\limits^{k_{n}}\Pr (J_{n}=k)\Pr \Big\{\mathop{\sup}\limits_{x\in[0,1]}|\sum_{j}^{k}\beta_{j}B_{j}|> M_{n}\Big\}+\mathop{\sum}\limits^{\infty}\limits_{k=k_{n}+1}\Pr (J_{n}=k)\\
&\lesssim \mathop{\sum}\limits_{k=1}\limits^{k_{n}}\Pr (J_{n}=k)k\exp \Big(-\frac{M_{n}^{2}}{2}\Big)+\mathop{\sum}\limits^{\infty}\limits_{k=k_{n}+1}\Pr (J_{n}=k)\\
&\lesssim k_{n}\exp \Big(-\frac{M_{n}^{2}}{2}\Big)+\mathop{\sum}\limits^{\infty}\limits_{k=k_{n}+1}\Pr (J_{n}=k)\\
&\lesssim e^{-n\epsilon^{2}}\end{split}
\end{equation*}
We derive the last $\lesssim$ through the facts that $k_{n}\exp \Big(-\frac{M_{n}^{2}}{2}\Big)\lesssim e^{-n\epsilon^{2}}$, and the assumption $ p_{n}^{k_{n}-1}(1-p_{n})=e^{-n\epsilon^{2}_{n}}$.

\subsection{Proof for theorem 2 }\label{proof No.3}
We denote $\kappa$ to be $\alpha$ or $\gamma$. By theorem 3.1 in \cite{AdaGau09}, there exists a Borel measurable subset $B_{n}^{(g)}$ of $C[0,1]^{d}$ such that 
\begin{equation}
\label{ConvCond}
\begin{split}
&\Pr(||W^{(g)A}-g_{0}||_{\infty}\leq \epsilon_{n})>e^{-n\epsilon_{n}^{2}};\\
&\Pr(W^{(g)A}\notin B_{n}^{(g)})\leq e^{-4n\epsilon_{n}^{2}};\\
&\log N(\epsilon_{n}, B_{n}^{(g)}, ||.||_{\infty})<K^{(g)}n\epsilon^{2}_{n},\\
\end{split}
\end{equation}\\
hold, for every sufficiently large n, and $\epsilon_{n}=n^{-\kappa/2(\kappa+d)}(\log n)^{(d+1)\kappa/(2\kappa+d)}$, $K^{(g)}$ is a sufficiently large constant. As stated in \cite{AdaGau09}, this power can be improved by using a slightly different prior for $A$. Then, the final rate of contraction will be improved, too, as which can be seen from the following proof.

We set $\Theta_{n}$ in the following way. Denote $\Theta^{(f)}_{n}=\{W^{A}\in B_{n}^{(f)}, and \ ||W^{A}||_{\infty}\leq N_{n}\}$. So, $\Theta^{(f)}_{n}$ increases to $B_{n}^{(f)}$ as n increases to infinity. As we assumed, $\{N_{n}\}$ is a sequence of real numbers increasing to infinity. We choose $N_{n}$ satisfying 
$$\Pr(W^{(f)A}\in B_{n}^{(f)})-\Pr(W^{A}\in \Theta_{n}^{(f)})\leq e^{-4n\epsilon_{n}^{2}}.$$
Then
$$\Pr(W^{A} \notin D^{(f)}_{n})\leq 2e^{-4n\epsilon_{n}^{2}}.$$
This can be achieved, since $\Pr(W^{(f)A}\in B_{n})$ goes to 1 and $e^{-4n\epsilon_{n}^{2}}$ goes to zero. Then we set $\Theta_{n}=B^{(\eta)}_{n}\times\Theta_{n}^{(f)}\subset C[0,1]^{d}\times C[0,1]^{d}$.

We start to verify all the conditions of general results on rates of posterior contraction. First, we bound the average hellinger distance entropy of the sieve of parameter 
space.
\begin{equation*}
\begin{split}
&\log N(\epsilon_{n}, \Theta_{n}, d_{n})\\
&\lesssim \log N(\epsilon_{n}, B^{(\eta)}_{n}, ||.||_{n})+\log N(\epsilon_{n},  \Theta^{(f)}_{n}, ||.||_{n})\\
&\leq \log N(\epsilon_{n}/e^{N_{n}}, B^{(\eta)}_{n}, ||.||_{\infty})+ \log N(\epsilon_{n}, B^{(f)}_{n}, ||.||_{\infty})\\
&\leq Kn\epsilon_{n}^{2}.
\end{split}
\end{equation*}
The first $\lesssim$ is from Lemma 1, the last $\leq$ is because of the third inequality of (\ref{ConvCond}). 

To estimate the prior positivity, we still use Lemma ~\ref{LemmaKL}. With the assumption that $||f||_{\infty}\leq N_{n}$, and $||f_{0}||_{\infty}\leq N_{n}$, for sufficiently large n, we can get
\begin{equation*}
\begin{split}
&\Pi_{n}(B_{n}((\eta_{0}, f_{0}),\epsilon_{n};2))\geq \Pi_{n}(||f-f_{0}||_{n}^{2}+||\eta-\eta_{0}||_{n}^{2}\leq e^{-4N_{n}}\epsilon_{n}^{2})\\
&\geq\Pi_{n}^{(f)}(f:||f-f_{0}||_{n}^{2}\leq \frac{e^{-4N_{n}}}{2}\epsilon_{n}^{2})\times \Pi_{n}^{(\eta)}(\eta:||\eta-\eta_{0}||_{n}^{2}\leq 
\frac{e^{-4N_{n}}}{2}\epsilon_{n}^{2})\\
&\geq \Pr(||W^{(f)A}-f_{0}||_{\infty}\leq \frac{e^{-2N_{n}}}{\sqrt{2}}\epsilon_{n})\times \Pr(||W^{(\eta)A}-\eta_{0}||\leq \frac{e^{-2N_{n}}}{\sqrt{2}}\epsilon)\\
&\geq e^{-2n\epsilon_{n}^{2}}.
\end{split}
\end{equation*}

Thus, for $\Theta_{n}\subset C[0,1]^{d}\times C[0,1]^{d}$ defined above, and $$\epsilon_{n}=\max\{ n^{-\alpha/(d+2\alpha)}(\log n)^{(d+1)\alpha/(2\alpha+d)}, n^{-\gamma/(d+2\gamma)}(\log n)^{(d+1)\gamma/(2\gamma+d)}\},$$
we have proved
\begin{equation*}
\begin{split}
&\log N(\epsilon_{n}, \Theta_{n}, d_{n})\leq 2Kn\epsilon_{n}^{2}\\
&\Pi_{n}(B_{n}((\eta_{0}, f_{0}),\epsilon;2))\geq  e^{-2n\epsilon_{n}^{2}}\\
&\Pi_{n}((f,\eta)\notin \Theta_{n})\leq 3e^{-4n\epsilon^{2}_{n}}.
\end{split}
\end{equation*}

The three assertions can be matched one-to-one with the assumption of general results on rates of posterior contraction (e.g. Theorem 4 in \cite{Ghosal00}), so the proof is completed.

The proof Corollary 2 is almost the same, except that the value of $\epsilon_{n}$ is given by Theorem 4.1 of \cite{ConGau08}.

\section{Discussion }\label{discussion}

In this paper, we investigated the posterior convergence rate for heteroscedastic nonparametric regression model with both mean function and variance function unknown and nonparametric. We considered both of the cases with random covariate  $\mathbf{x}$, and deterministic covariates. We also put the high-dimensional case in consideration. Though the rates we gave are not the minimax, they are only different with the optimal ones by a logarithmic factor. Besides, they are optimal for every regularity level. And we gave the minimax rate under the condition with known regularity level.

Whether the logarithmic factor of the posterior convergence rate is necessary for unknown regularity level is not known. To investigate this problem, other kinds of priors must be used, since as van der Vaart and  van Zanten have conjectured in \cite{AdaGau09}, the logarithmic factor is necessary with the rescaled Gaussian random field prior, and our current method used in the section of splines cannot  give the desired result, either.

\section{Appendix A. Proof of Lemma 1}\label{Appendix}
 By applying the inequalities $2-2ab\leq 2-2a+2-2b$, when $ a\leq 1$ and $\ b\leq 1$, together with $ 1-e^{-x}\leq x$ for $x\geq 0$, and $1-\frac{2x}{x^{2}+1}\leq  (2\log x)^{2}$ for all the x, we have 
 \begin{equation*}
 \begin{split}
&2-2exp(-\frac{(\eta_{1}(x)-\eta_{2}(x))^{2}}{4(V_{1}(x)+V_{2}(x))})\times \sqrt{\frac{2\sqrt{V_{1}(x)V_{2}(x)}}{V_{1}(x)+V_{2}(x)}}\\
&\leq 2 (1-\sqrt{\frac{2\sqrt{V_{1}(\mathbf{x})V_{2}(\mathbf{x})}}{V_{1}(\mathbf{x})+V_{2}(\mathbf{x})}})+2(1-exp\{-\frac{(\eta_{1}(\mathbf{x})-\eta_{2}(\mathbf{x}))^{2}}{4(V_{1}(x)+V_{2}(x))}\})\\
&\leq 2(\log (\frac{V_{1}(\mathbf{x})}{V_{2}(\mathbf{x})})^{2}+2 \frac{(\eta_{1}(\mathbf{x})-\eta_{2}(\mathbf{x}))^{2}}{4(V_{1}(x)+V_{2}(x))}.\\
\end{split}
\end{equation*}
Thus ,we have 
\begin{equation*}
 \begin{split}
d^{2}(\mathcal{P}_{\eta_{1},V_{1}},\mathcal{P}_{\eta_{2},V_{2}})&=\int\int(\mathcal{P}_{\eta_{1},V_{1}}^{\frac{1}{2}}-\mathcal{P}_{\eta_{2},V_{2}}^{\frac{1}{2}})^{2}\ dy\ dQ\\
&\leq 2 \int(\log (\frac{V_{1}(\mathbf{x})}{V_{2}(\mathbf{x})})^{2}+2 \frac{(\eta_{1}(\mathbf{x})-\eta_{2}(\mathbf{x}))^{2}}{4(V_{1}(x)+V_{2}(x))})\ dQ\\
\end{split}
\end{equation*}\\
held, which is followed by the result
$$\log N(3\epsilon, \Theta_{n}, d_{n}) \lesssim \log N(\epsilon/e^{N_{n}}, \Theta^{(\eta)}_{n},  ||.||_{n})+\log N(\epsilon, \Theta^{(f)}_{n},  ||.||_{n}).$$
provided $||V_{i}||>e^{-N_{n}}, i=1,2$.

\section{Appendix B. Proof of Lemma 2}\label{Appendix}
For the Kullback-Leibler divergence, we have,
\begin{equation*}
\begin{split}
K_{\mathbf{x}}(\mathcal{P}_{\eta_{1},V_{1}},\mathcal{P}_{\eta_{2},V_{2}})&=\frac{1}{2}\log\frac{V_{2}}{V_{1}}-\frac{1}{2}(1-\frac{V_{1}}{V_{2}})+\frac{1}{2}\frac{[\eta_{1}(\mathbf{x})-\eta_{2}(\mathbf{x})]^{2}}{V_{2}(\mathbf{x})}\\
&=\frac{1}{2}|(f_{2}(\mathbf{x})-f_{1}(\mathbf{x}))-\frac{1}{2}(1-e^{f_{1}(\mathbf{x})-f_{2}(\mathbf{x})})|+\frac{1}{2}\frac{[\eta_{1}(\mathbf{x})-\eta_{2}(\mathbf{x})]^{2}}{V(\mathbf{x})}.
\end{split}
\end{equation*}\\
We know that,
for $|z|\leq 2N,$
$$|z-1+e^{-z}|\leq |z|+|e^{-z}-1|\leq (e^{2N}+1)|z|;$$
when $z\geq 1$,$$(e^{2N}+1)|z|\leq (e^{2N}+1)z^{2},$$
when $z\leq 1$ $$|z-1+e^{-z}|\leq \sum_{n=2}^{\infty} z^{n}/2 \leq \frac{|z|^{2}/2}{1-|z|}\leq (e^{2N}+1)z^{2}.$$
Thus:
$$K(\mathcal{P}_{\eta_{1},V_{1}},\mathcal{P}_{\eta_{2},V_{2}})\leq (1+ e^{2N})(||\eta_{1}-\eta_{2}||_{n}^{2}+||(f_{1}-f_{2})||_{n}^{2}).$$

For the variance divergence, we have
$$Var_{x}(\mathcal{P}_{\eta_{1},V_{1}},\mathcal{P}_{\eta_{2},V_{2}})=2[-\frac{1}{2}+\frac{1}{2}\frac{V_{1}(x)}{V_{2}(x)}]^{2}+[\frac{V_{1}(x)}{V_{2}(x)}[\eta_{1}(x)-\eta_{2}(x)]]^{2}.$$\\
We can finish the proof with the inequality $|1-e^{z}|^{2}\leq (e^{2N})^{2}z^{2}$  for  $|z|\leq 2N$.



\bibliographystyle{elsarticle-num}
\bibliography{reference,papers.txt,books.txt}

\end{document}